\documentclass[12pt]{article}

\usepackage{amssymb,amsmath,exscale}

\usepackage{graphicx}

\usepackage{color}
\definecolor{marin}{rgb}   {0.,   0.3,   0.7} 
\definecolor{rouge}{rgb}   {0.8,   0.,   0.} 
\definecolor{sepia}{rgb}   {0.8,   0.5,   0.} 
\usepackage[colorlinks,citecolor=marin,linkcolor=rouge,
            bookmarksopen,
            bookmarksnumbered
           ]{hyperref}

\newtheorem{lemma}{Lemma}[section]

\newtheorem{theorem}[lemma]{Theorem}

\newtheorem{remark}[lemma]{Remark}
\newtheorem{example}[lemma]{Example}
\newtheorem{hypothesis}[lemma]{Hypothesis}
\newtheorem{notation}[lemma]{Notation}
\newtheorem{definition}[lemma]{Definition}
\newtheorem{conclusion}[lemma]{Conclusion}
\numberwithin{equation}{section}
\newcommand{\QED}{\mbox{}\hfill \raisebox{-0.2pt}{\rule{5.6pt}{6pt}\rule{0pt}{0pt}} 
          \medskip\par}             

\newenvironment{Proofof}[1]{\noindent
    \parindent=0pt\abovedisplayskip = 0.5\abovedisplayskip
    \belowdisplayskip=\abovedisplayskip{\bfseries Proof of #1. }}{\QED}

\newcommand{\Ac}{\mathcal{A}}

\newcommand{\dd}{\mathrm{d}}

\newcommand{\Hc}{\mathcal{H}}

\newcommand{\Ic}{\mathcal{I}}

\newcommand{\jb}{{\boldsymbol{j}}}

\newcommand{\N}{\mathbb{N}}
\newcommand{\Nc}{\mathcal{N}}

\newcommand{\Mc}{\mathcal{M}}
\newcommand{\Pc}{\mathcal{P}}

\newcommand{\R}{\mathbb{R}}

\newcommand{\C}{\mathbb{C}}

\newcommand{\T}{\mathbb{T}}
\newcommand{\Z}{\mathbb{Z}}
\newcommand{\Tc}{\mathcal{T}}

\newcommand{\Uc}{\mathcal{U}}
\newcommand{\Vc}{\mathcal{V}}

\newcommand{\Zc}{\mathcal{Z}}

\newcommand{\Norm}[2]{\|#1\|\left.\vphantom{T_{j_0}^0}\!\!\right._{#2}}         
\newcommand{\SNorm}[2]{|#1|\left.\vphantom{T_{j_0}^0}\!\!\right._{#2}}             

\title{Resonances in long time integration of semi linear Hamiltonian PDEs}        
\author{Erwan Faou{\small$\,^1$} and Beno\^it Gr\'ebert{\small$\,^2$}\\[4ex]
{\small$\,^1$}\, \small INRIA \& ENS Cachan Bretagne,  \\[-1ex]
\small Avenue Robert Schumann F-35170 Bruz, France. \\[-1ex]
\it \small email: \tt Erwan.Faou@inria.fr\\[2ex]
{\small$\,^2$}\, \small Laboratoire de Math\'ematiques Jean Leray,
Universit\'e de Nantes,\\[-1ex]
\small 2, rue de la Houssini\`ere
F-44322 Nantes cedex 3, France. \\[-1ex]
 \small and\\[-1ex]
 \small Instituto de Matematicas (Unidad Cuernavaca)\\[-1ex]
 \small Universidad Nacional Autónoma de México\\[-1ex]
 \small Av. Universidad, CP 62210, Cuernavaca, Mexico.\\[-1ex]
\small\it email: \tt benoit.grebert@univ-nantes.fr\\[2ex]
}       
\begin{document}
\maketitle

\abstract{
We consider a class of Hamiltonian PDEs that can be split into a linear unbounded operator and a regular non linear part, and we analyze their numerical discretizations by symplectic methods when the initial value is small in Sobolev norms. The goal of this work is twofold: First we show how standard approximation methods cannot in general avoid resonances issues, and we give numerical examples of pathological behavior for the midpoint rule and implicit-explicit integrators. Such phenomena can be avoided by suitable truncations of the linear unbounded operator combined with classical splitting methods. We then give a sharp bound for the cut-off depending on the time step. Using a new normal form result, we show the long time preservation of the actions for such schemes for all values of the time step, provided the initial continuous system does not exhibit resonant frequencies. \\[2ex]
}

\tableofcontents

\section{Introduction}

In this work, we consider Hamiltonian partial differential equations whose Hamiltonian functional $H  = H_0 + P$ can be split into a linear operator $H_0$ and a nonlinearity $P(u)$. Typical examples are given by Schr\"odinger equations and wave equations (see \cite{Bam07,Greb07,BG06,FGP1,FGP2} for similar framework). 

The understanding of the {\em geometric} numerical approximation of such equations over long time has recently known many progresses, see \cite{FGP1,FGP2,CHL08c,GL08b}. In essence these results show that under the hypothesis that the initial value is small and the physical system does not exhibit resonant frequencies, then the numerical solution associated with splitting methods induced by the decomposition of $H$ will remain small for very long time in high Sobolev norm, under a non resonance condition satisfied by the step size $h < h_0$. An analysis of this condition ensures the absence of resonances only if a strong CFL assumption is made between $h$ and the spectral parameter of space discretization. Outside the CFL regime, splitting methods may always exhibit resonances, in the sense where for some specific values of the stepsize, the conservation properties of the initial system are lost. 

To avoid these problems - which would require in practice the exact knowledge of the whole spectrum of the linear operator - it is traditionally admitted that the uses of {\em implicit} and symplectic schemes should help a lot. Indeed for {\em linear} PDEs for which  splitting methods induce resonances (see \cite{DF07}), it has been recently shown  (see \cite{DF08a}) that the use of implicit-explicit integrators based on a midpoint treatment for the unbounded part allows to avoid numerical instabilities. 

In this work, we show that in the nonlinear setting, it is in general hopeless to find a numerical method avoiding resonances unless the linear part of the equation is suitably truncated in the high frequencies. By numerical examples, we actually show how the midpoint rule applied to a Hamiltonian PDE either can yield numerical resonances and destroy the preservation properties of the initial system, or on the opposite break physical resonances and produce unexpected long time actions preservation.

We then show that if the linear operator is truncated in high frequencies, splitting methods yield numerical schemes preserving the physical conservation properties of the initial system without bringing extra numerical resonances between the stepsize $h$ and the frequencies of the system. In the case where this high frequencies cut-off corresponds to a CFL condition, we somewhat give a sharp bound for the CFL condition to avoid resonances extending the results in \cite{FGP1,FGP2} and \cite{CHL08c,GL08b}. 

Roughly speaking, the phenomenon can be explained as follows: for splitting schemes, resonances reflect the needle of a control of a small divisor of the form 
\begin{equation}
\label{Esmalldiv}
\exp(i h \Omega(\jb)) - 1
\end{equation}
where $\jb = (j_1,\ldots,j_r)$ is a multi-index, and where 
$$
\Omega(\jb)  =  \pm \omega_{j_1} \pm \cdots  \pm \omega_{j_r}
$$
is a (signed) sum of the frequencies $\omega_j$ of the linear operator $H_0$ indexed by a countable set of indices $j$ (the sign $\pm$ depends of the index $j_k$, see below for further details). 
The control of these small divisors up to an order $r$ ensures the control of Sobolev norms of the numerical solution on times of order $\varepsilon^{-r}$ for an initial value of order $\varepsilon$. 

Such a control can be made when $h \Omega(\jb) < \pi$ where \eqref{Esmalldiv} is essentially equivalent to $h\Omega(\jb)$ which can be controlled using classical estimates inherited from the generic assumption of absence of resonances in the physical system (see for instance \cite{Greb07,Bam07,BG06}). Such an assumption will be satisfied under a CFL condition bounding the frequencies in terms of the space discretization parameter, but it is worth noticing that it depends {\em a priori} on $r$ which is an arbitrary parameter (see \cite{FGP1,FGP2,GL08b}). 

Outside this CFL regime, resonances can always occur when $h \Omega(\jb) \simeq 2 \ell \pi$ for some $\ell \in \Z\notin\{0\}$. In \cite{FGP2}, it is shown however that such situation is exceptional, as it corresponds to very specific values of the time-step $h$. 

In the case of a general symplectic Runge-Kutta method or for implicit-explicit split-step method based on an implicit symplectic integrator for the linear part, the small divisors to be controlled take the form 
\begin{equation}
\label{EsmalldivRK}
\exp(i  \Psi(h,\jb)) - 1
\end{equation}
where 
$$
\Psi(h,\jb) = \pm \psi(h \omega_{j_1}) \pm \cdots  \pm \psi(h \omega_{j_r}) 
$$
and where $\psi$ is associated with the RK method. For the midpoint rule, the corresponding function is 
\begin{equation}
\label{Earctan}
\psi(x) = 2 \arctan\big(\frac{x}2\big). 
\end{equation}
Hence for these methods the control of the small divisors depends on the relation $\Psi(h,\jb) = 2\pi\ell$ which is a {\em nonlinear} relation between $h$ and the frequencies $\omega$. Of course a very restrictive CFL hypothesis stating that the $h\omega_{j_k}$ are {\em small} (e.g. of order $h$) ensures the control of these small divisors, but in general resonances occur. On numerical examples, it is indeed possible to exhibit $h$ such that $\Psi(h,\jb) = 0$ while $\Omega(\jb)  \neq 0$ (numerical resonances). Moreover, the somewhat opposite situation is also possible: the presence of physical resonances (i.e. $ \jb$ such that  $\Omega(\jb) = 0$) are {\em broken} by the action of $h$, yielding a $\Psi(h,\jb) \neq 0$ and unexpected regularity preservation. We refer to the section devoted to numerical examples for concrete illustrations. 

To avoid resonances in general situations, it seems very difficult to avoid the two following restrictions: 
\begin{itemize} 
\item $\psi(x)$ is linear, i.e. the numerical schemes is base on a splitting scheme or exponential integrator.  If this is not the case, we believe that the pathological behaviors described above can always be observed. 
\item There is a frequency cut-off in order to avoid situations where $h \Omega(\jb) = 2 \ell \pi$ with $\ell \neq 0$. Note that this does not exactly correspond to a CFL condition, as high frequencies are allowed to exist without restriction, but the action on the linear operator on these high frequencies is cancelled. 
\end{itemize}
In this work we give a sharp explicit bound for the cut-off in order to avoid numerical resonances and we prove a new normal form result for the corresponding splitting methods, yielding the preservation of the regularity of the numerical solution for a number of iterations of the form 
\begin{equation}
\label{Etime}
n < \varepsilon^{-1} h^{-N} + \varepsilon^{-r+1}
\end{equation}
where $r$ and $N$ are given parameters, and 
where $\varepsilon$ denotes the size of the initial value in Sobolev norm. In comparison with \cite{FGP1,FGP2,CHL08c,GL08b}, the difference is that the CFL condition is only imposed on the linear operator, and moreover this condition is independent of the approximation parameters $r$ and $N$, and is sharp in the sense where we can in general exhibit numerical resonances if it is violated.

As a conclusion, we would like to stress the following: the action of the midpoint rule (or of any symplectic RK methods) on the operator $H_0$ is equivalent to a smoothing in high-frequencies which amounts in some sense to a cut-off. What we show here is that the nonlinearity of the smoothing function (essentially based on the function $\arctan(x)$) introduces in general numerical instabilities. To avoid them, the user should better make the cut-off by himself! Note that this only  requires the a priori knowledge of bounds for the growth of the eigenvalues  of the linear operator $H_0$.

\section{Numerical examples}

\subsection{Actions-breaker midpoint}

We consider the cubic nonlinear Schr\"odinger equation
$$
i \partial_t u = - \Delta u + V \star u + |u|^2 u
$$
on the torus $\T^1 = \R/2\pi\Z$. We consider $V$ the function with Fourier coefficients
$
\hat V(k) = -\frac{2}{10 + 2k^2}. 
$
The frequencies of the operator $H_0$ in Fourier basis are hence 
\begin{equation}
\label{Efreq}
\omega_k = k^2 - \frac{2}{10 + 2k^2},\quad k \in \Z.
\end{equation} 
We first consider the implicit-explicit integrator defined as  
\begin{equation}
\label{Emidsplit}
u^{n+1} = R(-ihH_0) \circ \varphi_{P}^h(u^n)
\end{equation}
where $u^n$ are approximations of the exact solution $u$ at discrete times $t_n = nh$, and where
$$
R(z) = \frac{1 + z/2}{1 - z/2}
$$
is the stability function of the midpoint rule, and where 
$$
\varphi_P^h(u)(x) = \exp(-ih|u(x)|^2) u(x)
$$
is the exact solution of the nonlinear part. 
Note that the {\em unperturbed} integrator (i.e. without the nonlinear term) can be written in terms of Fourier coefficients (see \cite{DF08a})
\begin{equation}
\label{Eatan}
R(-ihH_0)_k = \frac{1 - ih \omega_k/2}{1 + ih \omega_k/2} = \exp(2i \arctan(h \omega_k/2)). 
\end{equation}
We consider as initial value the function 
\begin{equation}
\label{Einit}
u^0 = \frac{0.1}{2 - 2 \cos(x)} + 0.05( 2e^{2ix}-2 e^{5ix} +  3 e^{7ix}). 
\end{equation}
We use a collocation space discretization with $K = 100$ Fourier coefficients. 
In figure 1 and 2, we plot the evolution of the sum of the numerical actions $|\hat u^n_k|^2 + |\hat u^n_{-k}|^2$ for all $k \in \N$ in logarithmic scale with respect to the time. 
\begin{figure}[ht]
\label{figure1}
\begin{center}
\rotatebox{0}{\resizebox{!}{0.3\linewidth}{%
   \includegraphics{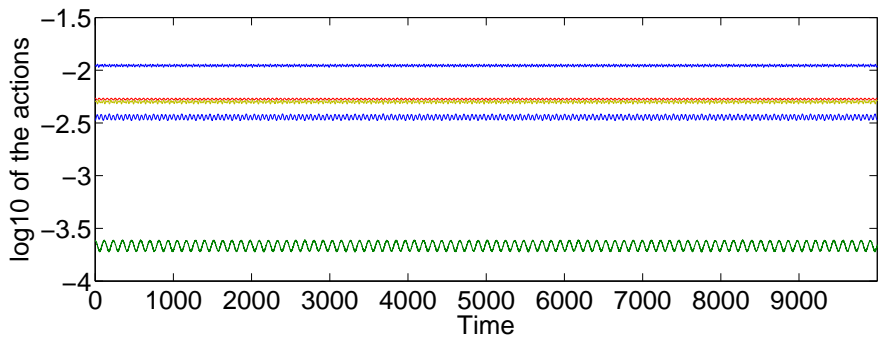} }}
\caption{\small{Plot of the actions for non-resonant step size (mid-split integrator).} }
\end{center}
\end{figure}
\begin{figure}[ht]
\label{figure2}
\begin{center}
\rotatebox{0}{\resizebox{!}{0.3\linewidth}{%
   \includegraphics{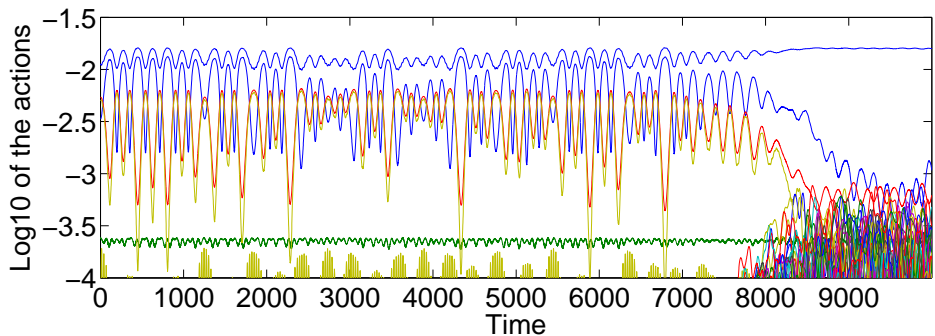}}}
\caption{\small{Plot of the actions for resonant step size  (mid-split integrator).} }
\end{center}
\end{figure}

In Figure 1, we use the stepsize $h = 0.13$, and we observe the long-time behavior of the actions. In Figure 2, we use a stepsize 
$h = 0.1278...$ such that 
\begin{equation}
\label{Ehreson}
\arctan(h\omega_2/2) + \arctan(h\omega_5/2) -  \arctan(h\omega_{-7}/2) = 0
\end{equation}
and we observe energy exchanges between the actions. This corresponds to the cancellation of the small divisor 
$$
\exp(2i \arctan(h\omega_2/2) + 2i\arctan(h\omega_5/2) -  2i\arctan(h\omega_{-7}/2) ) - 1
$$
which naturally appears in the search for a normal form for the numerical integrator \eqref{Emidsplit} (compare \cite{FGP1,FGP2}). 

We next consider the midpoint rule defined as 
\begin{equation}
\label{Emidpoint}
u^{n+1} = u^n -ihH_0 \big(\frac{u^{n+1} + u^{n}}{2}\big)   -ih g\big(\frac{u^{n+1} + u^{n}}{2}\big)
\end{equation}
where $g(u) = |u|^2u$. 
Equivalently, this scheme can be written
\begin{equation}
\label{Emidpoint2}
u^{n+1} = R(-ihH_0)\circ\big( u^n -\frac{2ih}{2 - ihH_0} g\big(\frac{u^{n+1} + u^{n}}{2}\big)\big).
\end{equation}
We use the same parameters as before. In Figure 3 we observe the preservation of the actions for $h = 0.13$ and in Figure 4 we see unexpected energy exchanges for the step size \eqref{Ehreson}. 
\begin{figure}[ht]
\label{figure3}
\begin{center}
\rotatebox{0}{\resizebox{!}{0.3\linewidth}{%
   \includegraphics{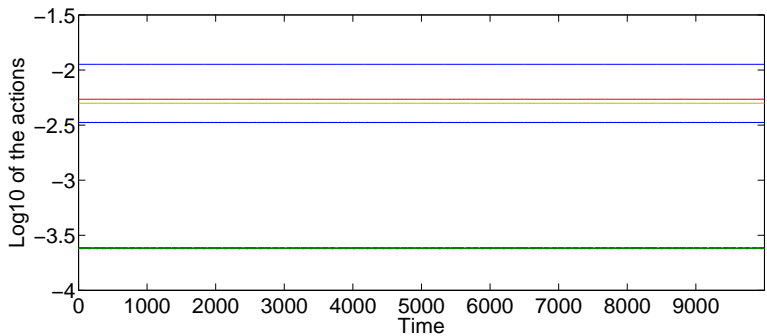} }}
\caption{\small{Plot of the actions for non-resonant step size (midpoint rule).} }
\end{center}
\end{figure}
\begin{figure}[ht]
\begin{center}
\rotatebox{0}{\resizebox{!}{0.3\linewidth}{%
   \includegraphics{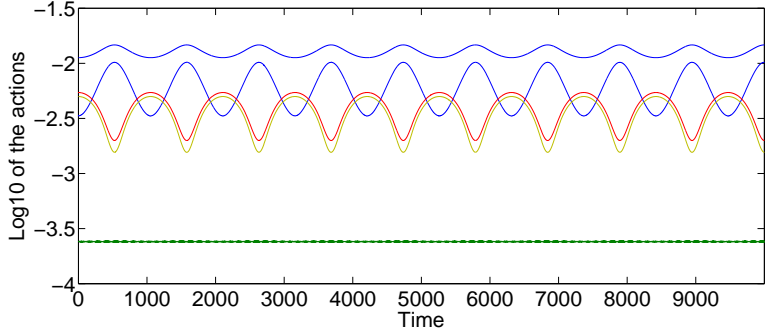}}}
\caption{\small{Plot of the actions for resonant step size (midpoint rule).} }
\end{center}
\end{figure}
Note that the numerical resonance effect is less pathological than for the mid-split integrator. This might be due to the regularization effect acting on the nonlinearity term (see \eqref{Emidpoint2}). However the preservation of the actions is lost anyway, reflecting the fact that the unperturbed linear integrator is the same as for the mid-split integrator, as well as the small divisors to be controlled.

\subsection{Resonance-breaker midpoint}

We consider now the same problem, but we set 
$$
\omega_2 = 10, \quad \omega_5 = 30 \quad \mbox{and}\quad \omega_{-7} = 40, 
$$
so that the (physical) resonance relation holds: 
$$
\omega_2 + \omega_5 - \omega_{-7} = 0. 
$$
The exact solution is plotted in figure 5. We observe energy exchanges between the actions. This solution is numerically computed using a splitting method under CFL condition with a CFL number close to 1.

\begin{figure}[ht]
\label{figure5}
\begin{center}
\rotatebox{0}{\resizebox{!}{0.3\linewidth}{%
   \includegraphics{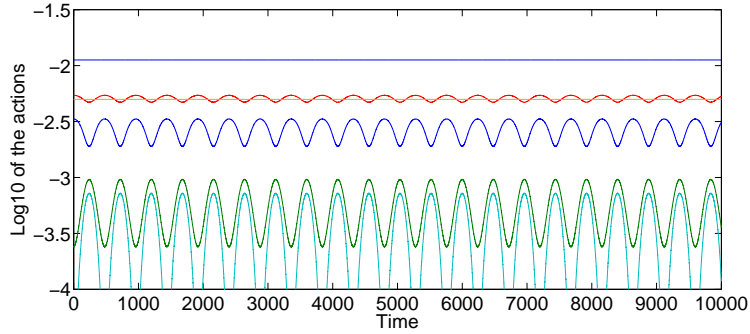} }}
\caption{\small{Exact solution with resonances.} }
\end{center}
\end{figure}

\begin{figure}[ht]
\label{figure6}
\begin{center}
\rotatebox{0}{\resizebox{!}{0.3\linewidth}{%
   \includegraphics{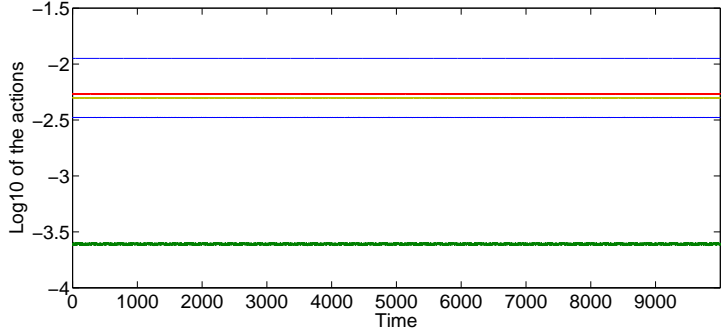}}}
\caption{\small{Undesired actions preservation by the mid-split integrator.} }
\end{center}
\end{figure}

In  Figure 6, we plot the evolution of the actions of the numerical solution given by the mid-split integrator \eqref{Emidsplit} with a step size $h = 0.0812...$ such that 
$$
2\arctan(h\omega_2/2) + 2\arctan(h\omega_5/2) -  2\arctan(h\omega_{-7}/2) = \frac12. 
$$
We observe the preservation of the actions over long time which reflects the control of the small divisor 
$$
|\exp( 2i\arctan(h\omega_2/2) + 2i\arctan(h\omega_5/2) -  2i\arctan(h\omega_{-7}/2)  ) - 1 | \simeq 0.485... 
$$ 
The same picture is plotted in Figure 7 for the midpoint rule \eqref{Emidpoint}. We thus see that these midpoint-like schemes are both unable to reproduce the qualitative behavior of the exact solution in Figure 5. 

\begin{figure}[ht]
\label{figure4}
\begin{center}
\rotatebox{0}{\resizebox{!}{0.3\linewidth}{%
   \includegraphics{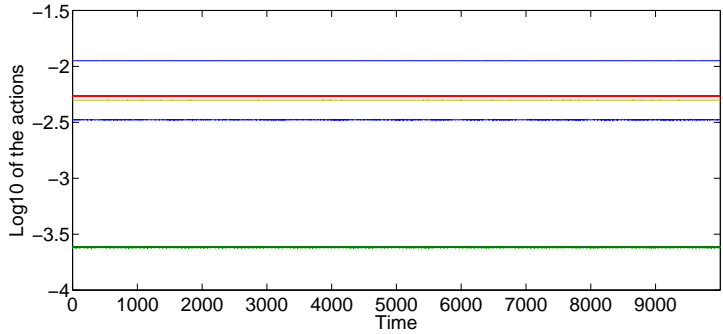}}}
\caption{\small{Undesired actions preservation by the midpoint rule.} }
\end{center}
\end{figure}

Note that the sets of stepsizes $h$ for which the solution looks like in Figure 6 or 7 for mid-split of midpoint rule integrators are very large. In fact, in this case the relation
\begin{equation}
\label{Esmd}
\Psi(h,\jb) = \pm 2 \arctan(h\omega_{j_1}/2) \pm \ldots \pm 2 \arctan(h\omega_{j_r}/2) = 0
\end{equation}
for a multi-index $\jb = (j_1,\ldots,j_r)$ 
 is exceptional, and not implied by a physical resonance 
 $$
 \Omega(\jb) = \pm \omega_{j_1} \pm \ldots \pm \omega_{j_r} = 0
 $$ 
 unless $h$ is very small so that the linear approximation of $\arctan(x)$ is valid. 

On the opposite, splitting methods will always catch these physical resonances.

\subsection{Conclusion and summary of the results}

The phenomena described above are somehow {\em generic} for any symplectic Runge Kutta method applied to Hamiltonian PDEs. Indeed for these methods, the unperturbed linear integrator can always  be written of the form \eqref{Eatan} (see for instance \cite{HLW}) and the control of the long-time behavior of the numerical solutions thus relies on the control of small divisors 
$$
\exp(i \Psi(h,\jb)) - 1
$$ 
where $\Psi(h,\jb)$ is of the form \eqref{Esmd}. This yields non linear relations between $h$ and the frequencies $\omega_k$ whose control is not implied by the physical small divisors $\Omega(\jb)$. 

To remedy this, a possibility is to use splitting methods (for which $\Psi(h,\jb) = h \Omega(\jb)$) with suitable truncation in high-frequencies. The goal of this paper is to give new results extending the existing ones in \cite{FGP1,FGP2,CHL08c,GL08b} by considering Hamiltonian PDEs where the operator $H_0$ is replaced by a truncated operator such that $h \omega_j \leq K$ form some $K$. This yields a truncated Hamiltonian $H_h = H_0^h + P$ where the nonlinear part is the same as for $H = H_0 + P$.  In this paper we show: 

\begin{itemize}
\item The preservation of the actions for the exact solution of the Hamiltonian PDE associated with $H_h$ under a generic non resonance condition on the frequencies of $H_0$. This preservation holds for times of the form $\varepsilon^{-1}h^{-N} + \varepsilon^{r +1}$ where $\varepsilon$ is the size of the solution in Sobolev norm. 
\item The persistence of this result for the numerical solution obtained by splitting methods induced by the decomposition of $H_h$ assuming an explicit bound for the cut-off level $K$. This preservation holds for all $h < h_0$ (without additional resonances) and for a number of iterations of the form \eqref{Etime}. The bound on $K$ turns out to be independent of $N$ and $r$ and sharp (see Thm. \ref{ET2} for a precise statement). 
\end{itemize}

We state the results in the Section 4, after recalling in Section 3 the mathematical framework already developed in  \cite{FGP2}. For ease of presentation,  we do not give the entire details of the proofs, as it uses very technical tools that might already be found in \cite{Greb07,FGP2} and \cite{BG06}. We prefer to give the key ingredients, and stress the specific changes in comparison with \cite{FGP2}.

\section{Abstract setting}
%
\subsection{Abstract Hamiltonian formalism}

This section essentially describes the formalism used in \cite{FGP1,FGP2}. 
We denote by $\Nc = \Z^d$ or $\N^d$ (depending on the concrete application) for some $d \geq 1$.  For $a = (a_1,\ldots,a_d) \in \Nc$, we set
$$
|a|^2 = \max\big(1,a_1^2 + \cdots + a_d^2\big). 
$$
We consider the set of variables $(\xi_a,\eta_b) \in \C^{\Nc} \times \C^{\Nc}$ equipped with the symplectic structure
\begin{equation}
\label{Esymp}
i \sum_{a \in \Nc} \dd \xi_a \wedge \dd \eta_a. 
\end{equation}
We define the set $\Zc = \Nc \times \{ \pm 1\}$. For $j = (a,\delta) \in \Zc$, we define $|j| = |a|$ and we denote by $\overline{j}$ the index $(a,-\delta)$. 

We will identify a couple $(\xi,\eta)\in \C^{\Nc} \times \C^{\Nc}$ with 
$(z_j)_{j \in \Zc} \in \C^{\Zc}$ via the formula
$$
j = (a,\delta) \in \Zc  \Longrightarrow 
\left\{
\begin{array}{rcll}
z_{j} &=& \xi_{a}& \mbox{if}\quad \delta = 1,\\[1ex]
z_j &=& \eta_a & \mbox{if}\quad \delta = - 1.
\end{array}
\right.
$$
By  a slight abuse of notation, we often write $z = (\xi,\eta)$ to denote such an element. 

For a given real number $s \geq 0$, we consider the Hilbert space $\Pc_s = \ell_s(\Zc,\C)$ made of elements $z \in \C^{\Zc}$ such that
$$
\Norm{z}{s}^2 := \sum_{j \in \Zc} |j|^{2s} |z_j|^2 < \infty,
$$
and equipped with the symplectic form \eqref{Esymp}. 

Let $\Uc$ be a an open set of $\Pc_s$. For a function $F$ of $\mathcal{C}^1(\Uc,\C)$, we define its gradient by 
$$
\nabla F(z) = \left( \frac{\partial F}{\partial z_j}\right)_{j \in \Zc}
$$
where by definition, we set for $j = (a,\delta) \in \Nc \times \{ \pm 1\}$, 
$$
 \frac{\partial F}{\partial z_j} =
  \left\{\begin{array}{rll}
 \displaystyle  \frac{\partial F}{\partial \xi_a} & \mbox{if}\quad\delta = 1,\\[2ex]
 \displaystyle \frac{\partial F}{\partial \eta_a} & \mbox{if}\quad\delta = - 1.
 \end{array}
 \right.
$$
Let $H(z)$ be a function defined on $\Uc$. If $H$ is smooth enough, we can associate with this function the Hamiltonian vector field $X_H(z)$ defined by
$$
X_H(z) = J \nabla H(z) 
$$
where $J$ is the symplectic operator on $\Pc_s$ induced  by the symplectic form \eqref{Esymp}.

For two functions $F$ and $G$, the Poisson Bracket is defined as
$$
\{F,G\} = \nabla F^T J \nabla G = i \sum_{a \in \Nc} \frac{\partial F}{\partial \eta_j}\frac{\partial G}{\partial \xi_j} -  \frac{\partial F}{\partial \xi_j}\frac{\partial G}{\partial \eta_j}.  
$$

We say that $z\in \Pc_s$ is {\em real} when $z_{\overline{j}} = \overline{z_j}$ for any $j\in \Zc$. In this case, $z=(\xi,\bar\xi)$ for some $\xi\in \C^{\Nc}$. Further we say that a Hamiltonian function $H$ is 
 {\em real } if $H(z)$ is real for all real $z$. 
 
\begin{definition}
Let $s \geq 0$, and let $\Uc$ be a neighborhood of the origin in $\Pc_s$. We denote by $\Hc^s(\Uc)$ the space of real Hamiltonian $H$ satisfying 
$$
H \in \mathcal{C}^{\infty}(\Uc,\C), \quad \mbox{and}\quad 
X_H \in \mathcal{C}^{\infty}(\Uc,\Pc_s). 
$$
\end{definition}

With a given function $H \in \Hc^s(\Uc)$, we associate the Hamiltonian system
$$
\dot z = J \nabla H(z)
$$
which can be written
\begin{equation}
\label{Eham2}
\left\{
\begin{array}{rcll}
\dot\xi_a &=& \displaystyle - i \frac{\partial H}{\partial \eta_a}(\xi,\eta) & a \in \Nc\\[2ex]
\dot\eta_a &=& \displaystyle i \frac{\partial H}{\partial \xi_a}(\xi,\eta)& a \in \Nc.  
\end{array}
\right.
\end{equation}
In this situation, we define the flow $\varphi_H^t(z)$ associated with the previous system (for times $t \geq 0$ depending on $z \in \Uc$). Note that if $z = (\xi,\bar \xi)$ and using the fact that $H$ is {\em real}, the flow $(\xi^t,\eta^t) = \varphi_H^t(z)$ satisfies for all time where it is defined the relation $\xi^t = \bar {\eta}^t$, where $\xi^t$ is solution of the equation 
\begin{equation}
\label{Eham1}
\dot\xi_a = - i \frac{\partial H}{\partial \eta_a}(\xi,\bar\xi), \quad a \in \Nc. 
\end{equation}
In this situation, introducing the real variables $p_a$ and $q_a$ such that
$$
\xi_a = \frac{1}{\sqrt{2}} (p_a + i q_a)\quad \mbox{and}\quad \bar{\xi}_a =  \frac{1}{\sqrt{2}} (p_a - i q_a),
$$
the system \eqref{Eham1} is equivalent to the system
$$
\left\{
\begin{array}{rcll}
\dot p_a &=& \displaystyle -  \frac{\partial H}{\partial q_a}(q,p) & a \in \Nc\\[2ex]
\dot q_a &=& \displaystyle  \frac{\partial H}{\partial p_a}(q,p),& 	a \in \Nc.  
\end{array}
\right.
$$
where $H(q,p) = H(\xi,\bar\xi)$. 

Note that the flow $\tau^t = \varphi_\chi^t$ of a real hamiltonian $\chi$ defines a symplectic map, i.e.  satisfies for all time $t$ and all point $z$ where it is defined
\begin{equation}
\label{Esympl}
(D_z  \tau^t)_z^T J (D_z\tau^t)_z = J
\end{equation}
where $D_z$ denotes the derivative with respect to the initial conditions. 

\subsection{Function spaces}

We describe now the hypothesis needed on the Hamiltonian $H$. 

Let $\ell \geq 3$ be a given integer.
We consider $ \jb = (j_1,\ldots,j_\ell) \in \Zc^{\ell}$, and we set for all $i = 1,\ldots p$,
$j_i = (a_i,\delta_i)$ where $a_i \in \Nc$ and $\delta_i \in \{\pm 1\}$. We define the moment $\Mc (\jb)$ of  the multi-index $\jb$ by
\begin{equation}
\label{EMcb}
\Mc(\jb) = a_{1} \delta_{1} + \cdots + a_\ell \delta_\ell.
\end{equation}
We then define the set of indices with zero moment
\begin{equation}
\label{EIr}
 \Ic_\ell =  \{  \jb = (j_1,\ldots,j_\ell) \in \Zc^{\ell}, \quad \mbox{with}\quad \Mc(\jb) = 0\}.
\end{equation}

 For $\jb = (j_1,\ldots,j_r) \in \Ic_r$, we define $\mu(\jb)$ as the third largest integer between $|j_1|,\ldots,|j_r|$. Then we set $S(\jb) = |j_{i_r}| - | j_{i_{r-1}}| + \mu(\jb)$ where $|j_{i_r}|$ and $|j_{i_{r-1}}|$  denote the largest and the second largest integer between $|j_1|,\ldots,|j_r|$. 

In the following, for $\jb = (j_1,\ldots,j_\ell) \in \Ic_\ell$, we use the notation 
$$
z_\jb = z_{j_1}\cdots z_{j_\ell}. 
$$
Moreover, for $\jb = (j_1,\ldots,j_\ell) \in \Ic_\ell$ with $j_i = (a_i,\delta_i) \in \Nc \times\{ \pm 1\}$ for $i = 1,\ldots,\ell$, we set
$$
\overline \jb = (\overline{j}_1,\ldots,\overline j_\ell)\quad\mbox{with}\quad \overline{j}_i = (a_i,-\delta_i), \quad i = 1,\ldots,\ell.
$$
We recall  the following definition from \cite{Greb07}. 
\begin{definition}
Let $k \geq 3$, $M > 0$ and $\nu \in [0,+\infty)$, and let
$$
Q(z) = \sum_{\ell = 3}^k \sum_{\jb \in \Ic_\ell} Q_{\jb} z_{\jb}.
$$
We say that $Q \in \Tc_k^{M,\nu}$ if there exist a constant $C$ depending on $M$ such that 
\begin{equation}
\label{Ereg}
\forall\, \ell = 3,\ldots,k,\quad \forall\, \jb \in \Ic_\ell,\quad |Q_{\jb}| \leq C \frac{\mu(\jb)^{M+\nu}}{S(\jb)^M}. 
\end{equation}
\end{definition}
Note that $Q$ is a real hamiltonian if and only if
\begin{equation}
\label{Ereal}
\forall\, \ell = 3, \ldots,k, \quad \forall\, \jb \in \Ic_\ell, \quad Q_\jb = \overline{Q}_{\overline\jb}. 
\end{equation}

We have that $\Tc^{M,\nu}_k \in \Hc^s$ for $s \geq \nu +1/2$ (see  \cite{Greb07}). The best constant in the inequality \eqref{Ereg} defines a norm $\SNorm{Q}{\Tc^{M,\nu}_k}$ for which $\Tc^{M,\nu}_k$ is a Banach space. 
We set
$$
T_k^{\infty,\nu} = \bigcap_{M \in \N} \Tc^{M,\nu}_k. 
$$


\begin{definition}
A function $P$ is in the class $\Tc$ if
\begin{itemize}
\item $P$ is a real hamiltonian and exhibits a zero of order at least 3 at the origin. 
\item There exists $s_0 \geq 0$ such that for any $s \geq s_0$, $P \in \Hc^s(\Uc)$ for some neighborhood $\Uc$  of the origin in $\Pc_s$. 
\item For all $k \geq 1$, there exists $\nu \geq 0$ such that the Taylor expansion of degree $k$ of $P$ around the origin belongs to $\Tc_k^{\infty,\nu}$. 
\end{itemize}
\end{definition}
With the previous notations, we consider in the following Hamiltonian functions of the form
\begin{equation}
\label{Edecomp}
H(z) = H_0(z) + P(z) = \sum_{a \in \Nc} \omega_a I_a(z)+ P(z),
\end{equation}
with $P \in \Tc$ and 
where for all $a\in \Nc$
$$
I_a(z) = \xi_a \eta_a
$$
are the {\em actions} associated with $a\in \Nc$. We assume that 
the frequencies $\omega_a \in \R$ satisfy \begin{equation}
\label{Eboundomega}
\forall\, a \in \Nc, \quad |\omega_a| \leq C |a|^m
\end{equation}
for some constants $C > 0$ and $m > 0$. 
The Hamiltonian system \eqref{Eham2} can hence be written
\begin{equation}
\label{Eham3}
\left\{
\begin{array}{rcll}
\dot\xi_a &=& \displaystyle - i \omega_a \xi_a - i \frac{\partial P}{\partial \eta_a}(\xi,\eta) & a \in \Nc\\[2ex]
\dot\eta_a &=& \displaystyle i \omega_a \eta_a + i \frac{\partial P}{\partial \xi_a}(\xi,\eta)& a \in \Nc.  
\end{array}
\right.
\end{equation}

\subsection{Non resonance condition}

Let  $\jb = (j_1,\ldots,j_r)\in \Ic_r$, and denote by $j_i = (a_i,\delta_i) \in \Nc \times \{\pm 1\}$ for $i = 1,\ldots,r$. We set 
\begin{equation}
\label{EbigOm}
\Omega(\jb) = 
\delta_1\omega_{a_1} + \cdots  + \delta_r\omega_{a_r}. 
\end{equation}
We say that $\jb = (j_1,\ldots,j_r) \in \Ic_r$ depends only of the actions and we write $\jb \in \Ac_r$ if $\jb = \overline{\jb}$. In this situation $r$ is even and we can write 
$$
\forall\, i = 1,\ldots r/2,\quad
j_{i} = (a_i,1), \quad\mbox{and}\quad j_{i + r/2} = (a_i,-1)
$$
for some $a_i \in \Nc$. 
Note that in this situation, 
$$
\begin{array}{rcl}
z_\jb = z_{j_1}\cdots z_{j_r} &=& \xi_{a_1}\eta_{a_1} \cdots \xi_{a_{r/2}} \eta_{a_{r/2}}\\[2ex]
&=& I_{a_1} \cdots I_{a_{r/2}}
\end{array}
$$
where for all $a \in \Nc$, 
$$
I_{a}(z) = \xi_a \eta_a
$$
denotes the action associated with the index $a$. Note that if $z$ satisfies the condition $z_{\overline{j}} = \overline{z_j}$ for all $j \in \Zc$, then we have $I_a(z) = |\xi_a|^2$. For odd $r$, $\Ac_r$ is the empty set.

We will assume now that the frequencies of the linear operator $H_0$ satisfy the following property:
\begin{hypothesis}\label{H1}
For all $r \in \N$, there exist constants $\gamma$ and $\alpha$ such that 
$\forall\,\jb = (j_1,\ldots,j_r) \notin  \Ac_r$, 
\begin{equation}
\label{nonres1}
|\Omega(\jb)| \geq \frac{\gamma}{\mu(\jb)^{\alpha}}. 
\end{equation}
\end{hypothesis}
where we recall that $\mu(\jb)$ denotes the third largest integer amongst $|j_1|, \ldots ,|j_r|$.

As explained in \cite{FGP1,FGP2,BG06,Greb07} a large set of  Hamiltonian PDEs can be written under the previous form and satisfy the non resonance condition \eqref{nonres1}. This includes the case of Schr\"odinger equations on a torus on dimension $d$ and wave equations with periodic boundary condition in dimension $1$.  
\subsection{Numerical integrator}

In the following $K$ is a fixed number and characterizes the high frequency cut-off level. As explained in the previous section the integrators we consider are such that the frequencies of the linear operator multiplied by $h$ are bounded by $K$. Let $\chi_K:\R_+ \to \R_+$ be the cut-off function 
$$
\chi_K(x) = 
\left\{\begin{array}{rl}
x & \mbox{if} \quad x \leq K\\[2ex]
0 & \mbox{if} \quad x > K
\end{array}
\right.
$$
We define the operator $ \psi_K(hH_0) := \exp(i \chi_K(h H_0))$ by the formula
$$
\forall\, j = (a,\delta) \in \Zc,\quad 
\psi_K(hH_0)_j = \exp(-i \delta \chi_K(h \omega_a)) = \left\{
\begin{array}{ll}
e^{-\delta ih\omega_a} & \mbox{if} \quad h\omega_a \leq K\\[2ex]
1 & \mbox{if} \quad h\omega > K
\end{array}
\right. 
$$
We consider split-step integrators of the form 
\begin{equation}
\label{Eint}
u^{n+1} = 
\psi_K(hH_0)
 \circ \varphi^h_P(u^n),
\end{equation}
where $\varphi^h_P$ is the exact flow of $P$.  

The numerical scheme \eqref{Eint} corresponds to an exact splitting method with time step $h$ applied to the truncated equation associated with the (infinite dimensional) Hamiltonian
\begin{equation}
\label{EHK}
H_h(z) = H_0^h(z) + P(z) := \sum_{a \in \Nc, |\omega_a| \leq Kh^{-1}} \omega_a \xi_a \eta_a + P(z), 
\end{equation}
where we denote by $H_0^h$ the truncated linear operator. 
The corresponding Hamiltonian system can be written
\begin{equation}
\label{Eham3K}
\left\{
\begin{array}{rcll}
\dot\xi_a &=& \displaystyle - i \frac1h\chi_K(h\omega_a) \xi_a - i \frac{\partial P}{\partial \eta_a}(\xi,\eta) & a \in \Nc\\[2ex]
\dot\eta_a &=& \displaystyle i \frac1h\chi_K(h\omega_a) \eta_a + i \frac{\partial P}{\partial \xi_a}(\xi,\eta)& a \in \Nc.  
\end{array}
\right.
\end{equation}

Note that the effect of $\chi_K$ is a cut-off in $\omega_a$ such that $\omega_a \leq K h^{-1}$. 

The following result is easily shown and given here without proof:
\begin{theorem}
\label{Tapprox}
Assume given two solutions $z(t)$ and $z_h(t)$ of \eqref{Eham3} and \eqref{Eham3K} respectively, such that $z_h(0) = z(0)$. Let $s,\sigma > 0$ and assume that for 
$0 \leq t \leq T$ we have $z(t)$ and $z_h(t)$ in $\Pc_{s+\sigma + \frac12}$. Then we have
$$
\forall\, t \in (0,T), \quad \Norm{z(t) - z_h(t)}{s} \leq C h^{\sigma}
$$
where $C$ depends on the Sobolev norms of $z(t)$ in $\Pc_{s+ \sigma+\frac12}$ for $t \in (0,T)$. 
\end{theorem}

\section{Statement of the result and applications}

We first give normal form results both for the truncated equation \eqref{Eham3K} and the numerical integrator \eqref{Eint}. We then give the dynamical consequences of these results for the long time behavior of the corresponding solutions.

\subsection{Normal form results}
\begin{theorem}
\label{ET1}
Assume that $P \in \Tc$, and that the non resonance condition \eqref{nonres1} is satisfied. Let $K$ and $r$ be given numbers. Then there exist $s_0$ and $h_0$ such that for all $s \geq s_0$ there exist $\Uc_s$ and $\Vc_s$ two neighborhoods of the origin in $\Pc_s$ such that for all $h \leq h_0$ there exists $\tau_s: \Vc_s \to \Uc_s$ a canonical transformation which is the restriction to $\Uc_s$ of  $\tau = \tau_{s_0}$ and which put the Hamiltonian $H_h$ of eqn. \eqref{EHK} under normal form 
$$
H_h \circ \tau = 
H_0^h + Z + R
$$
where $H_0^h$ is the Hamiltonian defined in \eqref{EHK} and where 
\begin{itemize}
\item[(i)] $Z$ is a real hamiltonian, polynomial of order $r$ in $z$ with terms that either depends only on the actions or contain at least three components $z_j$ with $|j| \geq Kh^{-1}$. As a consequence we have
\begin{equation}
\label{EZ1}
\forall\, z \in \Vc_s,\quad 
\sum_{a \in \Nc} |a|^{2s}|\{ I_a, Z\} (z)| \leq C h^s \big(\Norm{z}{s}^3 + \Norm{z}{s}^r\big)
\end{equation}
where $C$ depends on $r$, $s$ and $K$. 
\item[(ii)] $R \in \Hc^s(\Vc_s,\C)$ is a real hamiltonian such that for $z \in \Vc_s$, we have 
\begin{equation}
\label{ER1}
\Norm{X_R(z)}{s} \leq C \Norm{z}{s}^{r}
\end{equation}
where $C$ depends on $r$, $s$ and $K$. 
\item[(iii)] $\tau$ is close to the identity in the sense where for all $z \in \Vc_s$ we have 
\begin{equation}
\label{Etau1}
\Norm{\tau(z) - z}{s} \leq C \Norm{z}{s}^2
\end{equation}
and for all $z \in \Uc_s$
\begin{equation}
\label{Etau2}
\Norm{\tau^{-1}(z) - z}{s} \leq C \Norm{z}{s}^2
\end{equation}
where $C$ depends on $r$, $s$ and $K$.
\end{itemize}
\end{theorem}
As explained in the proof, this result is a mixed between the truncature systematically made in \cite{BG06} and the global result stated in \cite{Greb07}. The dependancy on $h$ in the estimates reflects the control of the non resonance conditions associated with the truncated linear operator appearing in $H_h$. 

\begin{theorem}
\label{ET2}
Assume that $P \in \Tc$, and that the non resonance condition \eqref{nonres1} is satisfied. Let $r$ be a given number. Let $K$ be such that 
\begin{equation}
\label{ECFL1}
K \leq \frac{\pi}{3},
\end{equation}
then there exist $s_0$ and $h_0$ such that for all $s \geq s_0$  there exist $\Uc_s$ and $\Vc_s$ two neighborhoods of the origin in $\Pc_s$ such that for all $h \leq h_0$ there exists a canonical transformation $\tau_s: \Vc_s \to \Uc_s$ which is the restriction to $\Uc_s$ of  $\tau = \tau_{s_0}$ such that 
\begin{equation}
\label{ENFsplit}
\tau^{-1} \circ \psi_K(hH_0)\circ \varphi_P^h \circ \tau = \psi_K(hH_0)\circ \theta
\end{equation}
where $\theta$ is the solution at time $\lambda = 1$ of a non-autonomous hamiltonian $h Z(\lambda) + R(\lambda)$ with
\begin{itemize}
\item[(i)] $Z(\lambda)$ a real hamiltonian depending smoothly on $\lambda$, polynomial of order $r$ in $z$ with terms that either depend only on the actions or contain at least three components $z_j$ with $|j| \geq (r-2)^{-1}Kh^{-1}$. 
\item[(ii)] $R(\lambda) \in \Hc^s(\Vc_s,\C)$ a real hamiltonian depending smoothly on $\lambda \in (0,1)$, and satisfying \eqref{ER1} uniformly in $\lambda$. 
\item[(iii)] $\tau$ close to the identity in the sense where it satisfies \eqref{Etau1} and \eqref{Etau2}. 
\end{itemize}
As a consequence, there exist a constant $C$ depending on $r$ and $s$ such that 
\begin{equation}
\label{Etheta}
 \forall\, z \in \Vc_s,\quad \sum_{a \in \Nc} |a|^{2s}|I_a(\theta(z)) - I_a(z) | \leq C \big(h^s \Norm{z}{s}^3 + \Norm{z}{s}^{r+1}\big). 
\end{equation}
\end{theorem}

\begin{remark}
As will appear clearly in the proof, the bound \eqref{ECFL1} can be refined to $\pi/2 - \delta$ with $\delta > 0$ for general situations. If the bound $\pi/2$ is not satisfied, we can construct a system such that numerical resonances between $h$ and the frequency vector $\omega$ appear. 
\end{remark}

\subsection{Dynamical consequences}

We now give the main outcome of the previous theorems: The first concerns the exact solution of \eqref{Eham3K} and the second the long time behavior of the numerical solution associated with splitting methods applied to this equation. 

\begin{theorem}
\label{ET3}
Assume that $P \in \Tc$ and $H_0$ satisfies the condition \eqref{nonres1}. Let $r,N \in \N^*$  be fixed. Then there exist constants $s_0$ and $h_0$ depending on $r$ and $N$ such that for all $s > s_0$, there exist a constant $\varepsilon_0$ depending on $r$, $N$, $K$ and $s$ 
such that 
the following holds: For all $\varepsilon < \varepsilon_0$, $h < h_0$ and for all $z^0\in \Pc_{s}$ real such that 
$\Norm{z^0}{s} < \varepsilon$, then the solution $z_h(t)$ of \eqref{Eham3K} with $z_h(0) = z^0$ satisfies 
\begin{equation}
\label{EresnormK}
\Norm{z_h(t)}{s} \leq 2 \varepsilon \quad \mbox{for}\quad t \leq  \frac{c}{\varepsilon^{r-1}} + \frac{c}{\varepsilon h^{N}},
\end{equation}
and  
\begin{equation}
\label{EresactK}
\sum_{a \in \Nc} |a|^{2s}| I_a(z_h(t)) - I_a(z_h(0))| \leq \varepsilon^{5/2}\quad \mbox{for}\quad t \leq \frac{c}{\varepsilon^{r-1}} + \frac{c}{\varepsilon h^{N}}. 
\end{equation}
for some constant $c$ depending on $s$, $r$ and $N$. 
\end{theorem}

Note that Eqn. \eqref{Eham3K} is a infinite dimensional PDE. The only frequency cut-off is made in the linear operator. The difference with the classical results \cite{BG06,Greb07} is the dependence of the cut-off parameter in the bound in time. For fully discretized systems  obtained by  pseudo spectral methods, the same result holds with constant independent of the spatial discretization parameter ({\em a priori} independent of $h$ and $K$). We do not give the proof here and refer to \cite{FGP1} for the description of fully discretized systems. 

\begin{theorem}
\label{ET4}
Assume that $P \in \Tc$ and $H_0$ satisfies the condition \eqref{nonres1}. Let $r,N \in \N^*$  be fixed, then there exist  constants $s_0$  and $h_0$ depending on $r$ and $N$ such that for all $s > s_0$, there exists a constant $\varepsilon_0$ depending on $r$, $N$ and $s$ 
such that 
the following holds: For all $\varepsilon < \varepsilon_0$, $h < h_0$ and for all $z^0\in \Pc_{s}$ real and 
$\Norm{z^0}{s} < \varepsilon$ 
if we define 
\begin{equation}
\label{Eseuil}
z^n = \big(\psi_K(h H_0) \circ \varphi^h_P \big)^n (z^0),
\end{equation}
where the frequency cut-off is such that
\begin{equation}
\label{ECFL}
K \leq \frac{\pi}{3}. 
\end{equation}
then
we have $z^n$ still real, 
and moreover
\begin{equation}
\label{Eresnorm}
\Norm{z^n}{s} \leq 2 \varepsilon \quad \mbox{for}\quad n \leq \frac{c}{\varepsilon^{r-1}} + \frac{c}{\varepsilon h^{N}},
\end{equation}
and  
\begin{equation}
\label{Eresact}
\sum_{a \in \Nc} |a|^{2s}| I_a(z^n) - I_a(z^0)| \leq \varepsilon^{5/2}\quad \mbox{for}\quad n \leq \frac{c}{\varepsilon^{r-1}} + \frac{c}{\varepsilon h^{N}}
\end{equation}
for some constant $c$ depending on $r$, $N$ and $s$. 
\end{theorem}

\begin{Proofof}{Theorem \ref{ET3}}
Let $y^0 = \tau^{-1}(z^0)$ which is well defined provided $\varepsilon_0$ is small enough so that $z^0 \in \Uc_s$. Using \eqref{Etau2} we have 
$$
\Norm{y^0 - z^0}{s} \leq C \varepsilon^2
$$
so that we can assume that $\Norm{y^0}{s} \leq \varepsilon/2$. 
We then define $y_h(t)$ the solution of the Hamiltonian system associated with the Hamiltonian $H_0^h + Z + R$ given in Theorem \eqref{ET1}. 
We set $\Nc(t) = \Norm{y_h(t)}{s}^2 = \sum_s j^{2s} I_j(y_h(t))$. We have 
$$
\frac{\dd \Nc}{\dd t} = \sum_s j^{2s} \{I_j(y_h(t)),Z(y_h(t)) + R(y_h(t)\} . 
$$  
Assume that  $s_0 > N$ and using \eqref{ER1} and \eqref{EZ1}, we see that as long as $y_h(t) \in \Vc_s$ we have 
$$
\left| \frac{\dd \Nc}{\dd t}\right| \leq C \big(h^{N} \Nc(t)^\frac{3}{2} + \Nc(t)^{\frac{r+1}{2}}\big)
$$
as we can always assume that $h_0 \leq 1$ and that $\Vc_s$ is contained is the ball of radius $1$ in $\Pc_s$. We know that $\Nc(0) \leq \varepsilon^{2}$ and we can assume that $\varepsilon_0$ is sufficiently small in such way that the ball in $\Pc_s$ centered at the origin  and of radius $3 \varepsilon/2$ is included in $\Vc_s$. Now as long as $\sqrt{\Nc(t)} \leq 3 \varepsilon/2$ we have
$$
\left| \frac{\dd \Nc}{\dd t}\right| \leq C (h^s \varepsilon^3 + \varepsilon^{r+1}) 
$$
for some constant $C$. This implies 
\begin{equation}
\label{Zug}
|\Nc(t) - \Nc(0)| \leq C t (h^N \varepsilon^3 + \varepsilon^{r+1}). 
\end{equation}
Hence there exists a constant $c$ such that as long as 
$$
t \leq c (h^{-N}\varepsilon^{-1} + \varepsilon^{r-1}) 
$$
we have $\sqrt{\Nc(t)} \leq 3\varepsilon/2$ and hence $y_h(t) \in \Vc_s$. But this implies that $z_h(t) = \tau(y_h(t))$. Using \eqref{Etau1} we then easily see that  \eqref{EresnormK} is satisfied. 
The proof of \eqref{EresactK} is similar (see \cite{BG06, Greb07}). 
\end{Proofof}

\begin{Proofof}{Theorem \ref{ET4}}
The proof follows the same lines as above but on a discrete level. Setting $y^0 = \tau^{-1}(z^0)$ we have $\Norm{y^0}{s} \leq \varepsilon/2$. Now if we define by induction
$$
y^{n+1} = \psi_K(hH_0)\circ \theta(y^{n})
$$
and if we define $\Nc_n = \Norm{y^n}{s}^{2}$
then estimate \eqref{Etheta} shows that 
$$
|\Nc_{n+1} - \Nc_n| \leq C (h^{N} \Nc_n^\frac{3}{2} + \Nc_n^{\frac{r+1}{2}}),
$$
provided $s_0 \geq N$, 
which is the discrete version of \eqref{Zug}. We then easily conclude upon using the same arguments. 
\end{Proofof}

\section{Proof of the normal form results}

\subsection{The continuous case}

We start now the proof of Theorem \ref{ET1}. 

The strategy follows lines of \cite{BG06,Greb07}: we search by induction a transformation eliminating the polynomial term of order $\ell$ until $\ell = r$. The transformation $\tau$ is then defined by the composition of all these transformations. At each step, the transformation is constructed as the flow at time $1$ of a Hamiltonian system with unknown Hamiltonian $\chi$. Hence we are lead  to solve recursively the homological equations
\begin{equation}
\label{Ehomc}
\{ H_0^h, \chi\} + Z = G
\end{equation}
where 
$$
G = \sum_{\jb \in \Ic_\ell} G_{\jb} z_\jb \in \Tc_\ell^{\infty,\nu_1}
$$ 
is a polynomial of order $\ell$ depending on the term constructed in the previous steps and satisfying estimates of the form \eqref{Ereg} for some $\nu_1 > 0$, and $Z$ an unknown term in normal form. 

Setting 
\begin{equation}
\label{Eomh}
\omega_a^h = \left\{ 
\begin{array}{ll}
\omega_a & \mbox{if}\quad \omega_a \leq K h^{-1}\\[2ex]
0 & \mbox{if}\quad \omega_a > K h^{-1}
\end{array}\right.
\end{equation}
we see that the equation \eqref{Ehomc} can be written in terms of the coefficients $\chi_\jb$, $Z_\jb$ and $G_{\jb}$ as
\begin{equation}
\label{Esyshom}
\Omega^h(\jb) \chi_\jb  + Z_\jb = G_\jb
\end{equation}
where $\Omega^h(\jb)$ is defined as \eqref{EbigOm} with respect to $\omega_a^h$. Following \cite{Greb07} in the case where no frequency cut-off is made (i.e. $\Omega^h(\jb) = \Omega(\jb)$), the condition \eqref{nonres1} ensures that the system \eqref{Esyshom} can be solved by putting the terms depending on the actions in $Z_\jb$ and solving the rest to construct $\chi_\jb$ by inverting $\Omega(\jb)$. It is then clear that $\chi_\jb$ belongs to some $\Tc_\ell^{\infty,\nu_2}$ for some $\nu_2$ which ensures the control of the regularity of the transformation. 

In our situation, it is clear that $\Omega^h(\jb)$ does not fulfill the condition \eqref{nonres1}: if for instance $(j_1,j_2,j_3)$ is a multi-index with all components greater than $Kh^{-1}$ in modulus, then $\Omega^h(\jb)$ is equal to zero. 

On the other hand, let $\jb = (j_1, \ldots, j_\ell) \in \Ic_\ell\backslash \Ac_\ell$ a multi-index with at most two indices greater than $Kh^{-1}$. We can always assume that these two big indices are $j_1$ and $j_2$ with $j_1 \geq j_2$ and that $\mu(\jb) = j_3$. 

\begin{itemize}
\item If both are greater than $Kh^{-1}$ we have in fact
$$
|\Omega^h(\jb)| = |\Omega(j_3,\ldots,j_\ell)| \geq \frac{\gamma}{\mu(\jb)^\alpha} 
$$
upon using \eqref{nonres1} unless $(j_3,\ldots,j_\ell) \in \Ac_{\ell-2}$. But in this last situation, the condition $\Mc(\jb) = 0$  implies that  $j_1 = \bar j_2$, i.e. $\jb \in \Ac_\ell$. 

\item If only one is greater than $Kh^{-1}$ then we have with similar notations
$$
|\Omega^h(\jb)| = |\Omega(j_2,\ldots,j_\ell)| \geq \frac{\gamma}{\mu(\jb)^\alpha} 
$$
thanks to \eqref{nonres1} unless $(j_3,\ldots,j_\ell) \in \Ac_{\ell-1}$. If $\ell$ is even this is impossible. If $\ell$ is odd, then the zero moment condition $\Mc(\jb) = 0$ implies that $j_1 = 0$ which is a contradiction. 

\end{itemize}

This shows that \eqref{nonres1} holds for $\Omega^h(\jb)$ except for $\jb \in \Ac_\ell$ or for $\jb$ such that at least three indices are greater that $Kh^{-1}$ in modulus. Hence we solve \eqref{Esyshom} by defining  $Z_\jb=G_\jb$ and $\chi_\jb=0$ when $\jb \in \Ac_\ell$ or when $\jb$ contains at least three indices are greater that $Kh^{-1}$ in modulus; while $Z_\jb=0$ and $\chi_\jb  = \Omega^h(\jb)^{-1} G_\jb$ in the other cases, i.e. when $\jb \notin \Ac_\ell$ and  $\jb$ contains  at most two indices are greater that $Kh^{-1}$ in modulus.

Now the condition \eqref{nonres1} ensures that $\chi_{\jb} \in \Tc_\ell^{\infty,\nu_2}$ for some $\nu_2$. The conclusion now follows \cite{Greb07} except the derivation of equation \eqref{EZ1} which is a consequence of Lemma 4.11 in \cite{BG06}.

\subsection{Splitting methods}

We prove now Theorem \ref{ET2}. 

We follow now the methodology developed in \cite{FGP1,FGP2}. We embed the splitting method $\Psi_K(hH_0) \circ \varphi_P^h$ into the family 
$$
(0,1) \ni \lambda \mapsto \psi_K(hH_0) \circ \varphi^\lambda_{hP}
$$
and we seek $\tau^\lambda$ as a transformation $\tau = \varphi_{\chi(\lambda)}^\lambda$ associated with a non-autonomous real hamiltonian $\chi(\lambda)$ depending smoothly on $\lambda$ and such that for all $\lambda \in (0,1)$
\begin{equation}
\label{eq:flots}
\forall\, \lambda \in [0,1]\quad \psi_K(hH_0) \circ \varphi_{hP}^\lambda \circ \varphi_{\chi(\lambda)}^\lambda =  \varphi_{\chi(\lambda)}^\lambda \circ\psi_K(hH_0) \circ \varphi_{hZ(\lambda)}^\lambda
\end{equation}
where $Z(\lambda)$ is a Hamiltonian in normal form in the sense of Theorem \ref{ET2} and $R(\lambda)$ a real Hamiltonian possessing a zero of order $r+1$. 
Deriving this expression in $\lambda$, we find the equation (compare eqn. (5.18) in \cite{FGP2}): 
\begin{equation}
\label{eq:tg2}
\forall\, \lambda \in [0,1]\quad 
  \chi(\lambda) \circ \psi_K(hH_0) - \chi(\lambda) \circ \varphi_{hP}^{-\lambda}= hP - (h Z(\lambda) + R(\lambda)) \circ \varphi^{-\lambda}_{\chi(\lambda)} . 
\end{equation}
As in \cite{FGP2}, we see that the solution of this equation relies on the solvability of a discrete Homological equation of the form 
$$
\forall\, \jb \in \Ic_\ell, \quad 
(e^{ih \Omega^h(\jb)} - 1)\chi_\jb + hZ_\jb = hG_\jb
$$
where now $\Omega^h(\jb)$ is defined as in \eqref{EbigOm} with respect to  $\omega_a^h$ defined in \eqref{Eomh}. 

Assume that $G_\jb \in \Tc_\ell^{\infty,\nu_1}$ for some $\nu_1$. 

Let $\jb \in \Ic_\ell\notin\Ac_\ell$. As before, we assume that $|j_1| \geq |j_2| \geq |j_3|\geq\cdots$. 
We recall that for $|x|\leq \pi$, we have 
$$
|e^{ix} - 1| \geq \frac2\pi |x|. 
$$
Assume that $\mu(\jb) = j_3 \leq (\ell-2)^{-1} K h^{-1}$. Then we have 
$$
|h \Omega(j_3,\ldots,j_\ell) | \leq K \leq \frac{\pi}{3}.
$$
Now we have three possibilities: 
\begin{itemize}
\item $|j_2| > Kh^{-1}$. In this situation we have $\omega^h_{j_1} = \omega^h_{j_2} = 0$. 
$$
|\Omega^h(\jb)| = |h \Omega(j_3,\ldots,j_\ell)| \leq \frac{\pi}{3}
$$
and hence
$$
|e^{ih \Omega^h(\jb)} - 1| \geq \frac{2}{\pi} h \Omega^h(j_3,\ldots,j_\ell) \geq \frac{h\gamma^*}{\mu(\jb)^\nu} 
$$
for some constant $\gamma^*$ using \eqref{nonres1} and unless $(j_3,\ldots,j_\ell) \in \Ac_{\ell-2}$. But in this last case, thanks to the  the zero moment condition, $j_1=\bar j_2$ and thus $\jb \in  \Ac_{\ell}$. 
\item $|j_1| > Kh^{-1}$ and $|j_2| \leq Kh^{-1}$.
Now we have 
$$
|\Omega^h(\jb)| = |h \Omega(j_2,\ldots,j_\ell)| \leq 2K \leq \frac{2\pi}{3}
$$  
and hence 
$$
|e^{ih \Omega^h(\jb)} - 1| \geq \frac{2}{\pi} h \Omega^h(j_2,\ldots,j_\ell) \geq \frac{h\gamma^*}{\mu(\jb)^\nu} 
$$
unless $(j_3,\ldots,j_\ell) \in \Ac_{\ell-2}$ which is impossible because the zero moment condition would be violated. 
\item $|j_1| \leq Kh^{-1}$. In this case we have 
$$
|\Omega^h(\jb)| = |\Omega(\jb)| \leq 3K \leq \pi
$$
and hence we have 
$$
|e^{ih \Omega^h(\jb)} - 1| \geq \frac{2}{\pi} h \Omega(\jb) \geq \frac{h\gamma}{\mu(\jb)^\nu} 
$$
using \eqref{nonres1}. 
\end{itemize}

So far we have proven the following: For all $\jb$ such that $\mu(\jb) \leq (\ell-2)^{-1} K h^{-1}$ and $\jb \notin  \Ac_{\ell}$, we have 
$$
|e^{ih \Omega^h(\jb)} - 1| \geq \frac{h\gamma^*}{\mu(\jb)^\nu}
$$
for some constant $\gamma^*$ depending on the constant $\gamma$ in \eqref{nonres1} (for different $\ell$). 

Now we see that we can solve the homological equation for those indices, and put the remainder terms in $Z$, which entails into the definition of normal form terms in Theorem \eqref{ET3}. It is then clear that $\chi(\lambda) \in \Tc_\ell^{\infty,\nu_2}$ for some $\nu_2$, and the rest of the proof can be adapted from \cite{Greb07}. 

Finally,  estimate \eqref{Etheta} is a direct consequence of Lemma 4.11 in \cite{BG06}. This concludes the proof. 
\section*{Acknowledgement}

The authors are glad to thank the Mathematical Institute of Cuernavaca (UNAM Mexico) where this work was initiated.

\end{document}